\documentclass{amsart}
\usepackage{bbm}\newcommand{\dontprint}[1]{\relax}
\newcommand{\R}{\mathbbm{R}}
\newcommand{\cA}{\mathcal{A}_d}

\newcommand{\del}[2]{\frac{\partial{#1}}{\partial{#2}}}
\newcommand{\Cdot}{{\displaystyle\cdot}}
\newcommand{\Ma}{M^{\mathrm{aff}}}
\newcommand{\Mco}{M^{\mathrm{coor}}}
\newtheorem%
{thm}{Theorem}[section]
\newtheorem%
{proposition}[thm]{Proposition}
\newtheorem%
{lemma}[thm]{Lemma}
\newtheorem%
{lemmadef}[thm]{Lemma-Definition}
\newtheorem%
{corollary}[thm]{Corollary}
\title[Fedosov connections and deformation quantization]
{Fedosov connections on jet bundles and deformation
quantization}
\author{A. S. Cattaneo${}^*$}
\author{G. Felder${}^{**}$}
\author{L. Tomassini${}^{**}$}
\address{A. S. C.: Department of Mathematics, Harvard University,
Cambridge,
MA~02138, and
Institut f\"ur Mathematik, Universit\"at Z\"urich,  
CH-8057 Z\"urich, Switzerland}
\email{asc@math.unizh.ch}
\address{G. F.: D-MATH, ETH-Zentrum, CH-8092 Z\"urich , Switzerland}
\email{felder@math.ethz.ch}
\address{L. T.: D-MATH, ETH-Zentrum, CH-8092 Z\"urich , Switzerland}
\email{lorenzo@math.ethz.ch}
\thanks{${}^*$
Work
partially supported by SNF Grant
No.~20-63821.00}\thanks{${}^{**}$Work
partially supported by SNF Grant
No.~21-65213.01}

\begin{document}
\begin{abstract} We review our construction of star-products on
Poisson manifolds and discuss some examples. In particular,
we work out the relation with Fedosov's original construction
in the symplectic case.
\end{abstract}

\maketitle
\section{Introduction}\label{s-1}
In \cite{K} Kontsevich solved the problem of describing 
associative deformations 
of the algebra of functions on a manifold. In particular,
he gave an explicit formula for the associative 
deformations 
(``star-products'' \cite{BFFLS})
$f\star g=fg+\epsilon B_1(f,g)+\epsilon^2B_2(f,g)+\cdots$
of the product of functions on an open subset $M$ of $\R^d$,
with bidifferential operators $B_j$ and so that $f\star 1=1\star f=f$.
Here $\epsilon$ is a formal parameter ($\epsilon=
\mathrm{i}\hbar/2$ in the notation of physics).
He also described a (non-explicit) way to 
obtain star-products on
arbitrary manifolds, as a corollary of his formality theorem, see \cite{K} and, for more details,
Appendix A.3 of \cite{K2}.

Here we review a more explicit version \cite{CFT} 
of the construction
of star-products on arbitrary manifolds. 
For a more physics-oriented review and the relation
to quantum field theory, see \cite{CF2}.
Our approach is based on the
notion of Fedosov connections on jet bundles.
A special case of this notion was introduced by Fedosov \cite{F} 
to construct star-products such that the Poisson bracket
defined by $f\star g-g\star f=2\epsilon\{f,g\}+O(\epsilon^2)$ comes
from a symplectic structure.

In general, one is given a Poisson manifold $(M,\alpha)$,
which is a manifold $M$ with a bivector field $\alpha\in
\Gamma(M,\wedge^2TM)$ so that the bracket $\{f,g\}=\alpha(df\otimes dg)$ obeys the Jacobi identity and asks for star-products
as above so that $\frac12(B_1(f,g)-B_1(g,f))=\{f,g\}$.

To construct such products using Kontsevich's local formula,
we first describe the local data given by the Poisson
structure. Let $E_0$ be the bundle of (infinite)
jets of smooth
functions on $M$. The fiber over $x\in M$ consists of
jets of functions at $x$, namely
equivalence classes of smooth functions defined on
some open neighborhood of $x$, where two functions are
considered to be equivalent iff they have the same Taylor
expansion at $x$ (with respect to any choice of
local coordinates around $x$). Similarly, one considers
jets of multivector fields and multidifferential operators.
To any function $f\in C^\infty(M)$ there
corresponds a section of $E_0$, the jet of $f$: 
its value at $x$ is the jet of $f$ at $x$. 
The bundle $E_0$ comes with a canonical flat connection
$D_0:\Gamma(M,E_0)\mapsto \Omega^1(M,E_0)$, which is the
unique connection so that its 
horizontal sections are precisely the
jets of globally defined smooth functions. Moreover
the jets of $\alpha$ define a Poisson bracket on
each fiber and we have the Leibniz rules
\[
D_0(fg)=D_0(f)g+fD_0(g),\qquad
D_0\{f,g\}=\{D_0(f),g\}+\{f,D_0(g)\},
\]
for any sections $f,g\in\Gamma(M,E_0)$.
Thus $E_0$ is a bundle of Poisson algebras with flat 
connection. The Leibniz rules imply that 
we have an isomorphism of Poisson algebras
\[
\iota:C^\infty(M)\to H^0(E_0,D_0)=\mathrm{Ker}(D_0),
\]
from the algebra of smooth functions onto the algebra of
horizontal sections.

The idea is now to ``quantize'' $(E_0,D_0)$: we construct
a bundle $E=E_0\otimes\R[[\epsilon]]$ 
of associative $\R[[\epsilon]]$-algebras
using Kontsevich's local formula in each fiber. Then
we deform the connection $D_0$ to a flat connection
$\bar D=D_0+O(\epsilon)$ on $E$ obeying the Leibniz
rule
\[
\bar D(f\star g)=\bar D(f)\star g+f\star \bar D(g),
\qquad f,g\in \Gamma(M,E).
\]
Thus the product induces a 
product on the 
space of horizontal sections $H^0(E,\bar D)$.
Finally we construct a {\em quantization map}, an
isomorphism of $\R[[\epsilon]]$-modules
$\rho:H^0(E_0,D_0)[[\epsilon]]\to H^0(E,\bar D)$
given by differential operators and such that
 $\rho(1)=1$.
Then $f\star_M g=\rho^{-1}(\rho(f)\star\rho(g))$ is a 
star-product on $M$.

\section{Formality theorem for $\R^d$}\label{s-2}
Let $\cA=\R[[y^1,\dots,y^d]]$ be the algebra 
of formal power series in $d$ indeterminates. A (formal)
Poisson bracket on $\cA$ is a Lie algebra
structure on $\cA$ of the form
\[
\{f,g\}=\sum_{i,j=1}^d\alpha^{ij}\del f{y^i}\del g{y^j},
\qquad f,g\in\cA,
\]
for some $\alpha^{ij}\in\R[[y^1,\dots,y^d]]$.
M. Kontsevich constructed in \cite{K} an $\R[[\epsilon]]$-bilinear
associative product $\star_\alpha$ on $\cA[[\epsilon]]$ which is
a deformation of the product on $\cA$: for $f,g\in\cA$,
the product has the form
\[
f\star_\alpha g=fg+\sum_{n=1}^\infty
\frac{\epsilon^n}{n!}U_n(\alpha,\dots,\alpha)(f\otimes g),
\]
where
$U_n(\alpha,\dots,\alpha)\in\mathrm{Hom}_\R(\cA\otimes\cA,\cA)$
is a certain bidifferential operator whose coefficients
are homogeneous polynomials of degree $n$ 
in the partial derivatives of the coordinates $\alpha^{ij}$
of the Poisson bivector field $\alpha$.
For example, $U_1(\alpha)(f\otimes g)=\sum_{i,j=1}^d
\alpha^{ij}\partial_if
\partial_jg=\{f,g\}$.

More generally, $U_n(\alpha_1,\dots,\alpha_n)$ is
defined not only for Poisson bivector fields but also
for multivector fields (skew-symmetric contravariant tensor
fields) $\alpha_i$. If $\alpha_i$ has rank $m_i$, $i=1,\dots,n$,
then $U_n(\alpha_1,\dots,\alpha_n)\in \mathrm{Hom}_{\R}
(\cA^{\otimes m},\cA)$, where $m=\sum m_i-2n+2$, is a
multidifferential operator. The coefficients of these
multidifferential operators are polynomials in
the partial derivatives of the coordinates
of $\alpha_1,\dots,\alpha_n$, linear in each $\alpha_i$.
They are given in terms of integrals of differential
forms over configuration spaces of $n$-points in the upper
half plane, see \cite{K}, and can be interpreted as Feynman amplitudes
for a topological string theory \cite{CF}.

The associativity of the product is a special case of
Kontsevich's {\em formality theorem}, which is a sequence
of quadratic relations for the operators $U_n$. For
our purpose we need a further set of special cases, namely
all cases where $\alpha_i$ is either a given Poisson
bivector field, a vector field or a function. Leaving
functions apart for the time being, the relations may be
expressed as follows: introduce the generating series
\begin{eqnarray}
P(\alpha)&=&\sum_{n=0}^\infty 
\frac{\epsilon^n}{n!}U_{n}(\alpha,\dots,\alpha)\in
\mathrm{Hom}_\R(\cA\otimes\cA,\cA),\label{e-P}\\
A(\xi,\alpha)&=&\sum_{n=0}^\infty 
\frac{\epsilon^{n}}{n!}U_{n+1}(\xi,\alpha,\dots,\alpha)\in
\mathrm{Hom}_\R(\cA,\cA),\label{e-A}\\
F(\xi,\eta,\alpha)&=&\sum_{n=0}^\infty 
\frac{\epsilon^n}{n!}
U_{n+2}(\xi,\eta,\alpha,\dots,\alpha)\in
\cA.\label{e-F}
\end{eqnarray}

Let $W_d=\oplus_{1}^d\cA\partial/\partial y^i$ 
be the Lie algebra of formal vector fields on $\R^d$. This
Lie algebra 
acts on $\cA$ and thus on the Hochschild cochain complex
$\mathrm{Hoch}^\Cdot(\cA,\cA)=\oplus_{n=0}^\infty
\mathrm{Hom}_\R(\cA^{\otimes n},\cA)$. Let 
\[
C^{\Cdot}_{\mathrm{Lie}}
(W_d,\mathrm{Hoch}^\Cdot(\cA,\cA))
=
\mathrm{Hom}_\R
(\wedge^{\Cdot} W_d,\mathrm{Hoch}^\Cdot(\cA,\cA))
\]
be the corresponding Lie algebra complex, with differential
$\delta$. The maps $\xi_1\wedge\cdots\wedge\xi_k\mapsto
U_{n+k}(\xi_1,\dots,\xi_k,\alpha,\dots,\alpha)$, for
$k=0,1,2$ can be then considered as elements of
this complex.
The formality theorem for these objects may be written
as
\noindent
\begin{enumerate}
\item[(i)] $P(\alpha)\circ(P(\alpha)\otimes \mathrm{Id})=
P(\alpha)\circ(\mathrm{Id}\otimes P(\alpha))$.
\item[(ii)]$
P(\alpha)\circ(A(\xi,\alpha)\otimes \mathrm{Id}+\mathrm{Id}\otimes A(\xi,\alpha))
-A(\xi,\alpha)\circ P(\alpha)
=\delta P(\xi,\alpha)$.
\item[(iii)]$
P(\alpha)\circ(F(\xi,\eta,\alpha)\otimes\mathrm{Id}-\mathrm{Id}\otimes F(\xi,\eta,\alpha))
-
A(\xi,\alpha)\circ A(\eta,\alpha)+A(\eta,\alpha)\circ A(\xi,\alpha)
=\delta A(\xi,\eta,\alpha)$.
\item[(iv)]
$
-A(\xi,\alpha)\circ F(\eta,\zeta,\alpha)
-A(\eta,\alpha)\circ F(\zeta,\xi,\alpha)
-A(\zeta,\alpha)\circ F(\xi,\eta,\alpha)
\!=\!\delta F(\xi,\eta,\zeta,\alpha)$.
\end{enumerate}
The first equation is the associativity of the product.
The second equation says that changing coordinates
leads to an equivalent product. 
Indeed $\delta P(\xi,\alpha)=\frac d{dt}P(\phi_{t*}\alpha)|_{t=0}$ is the
infinitesimal variation of the product under
a coordinate transformation given by the flow
$\phi_t$ of $\xi$.
Before explaining the meaning of the remaining equations,
which will appear in Prop.\ \ref{p-34} below,
we need to discuss the action of the Lie algebra $gl_d
\subset W_d$ of linear vector fields and the lowest order
terms of $P$, $A$, $F$.
The action of $gl_d$ on $P$, $A$, $F$ is described by
two properties:
\begin{enumerate}
\item[A.] The operators $U_n(\alpha_1,\dots,\alpha_n)$
are $gl_d$ equivariant, in the sense that
\[
g^*\circ
U_n(g_*\alpha_1,\dots,g_*\alpha_n)=U(\alpha_1,\dots,\alpha_n)
\circ(g^*)^{\otimes m},\]
where $g^*f(y)=f(gy)$ and
$g_*\alpha_i(y)=(g\otimes\cdots\otimes g)\alpha_i(g^{-1}y)$,
if $g\in\mathrm{GL}(d,\R)$, $f\in \cA$, and $\alpha_i$ is
a formal multivector field. Therefore $P,A,F$ are
 $gl_d$-equivariant.
\item [B.] For any $\xi\in gl_d$, $\eta\in W_d$,
$
A(\xi,\alpha)=\xi$, viewed as a first order differential
operator and $F(\xi,\eta,\alpha)=0$. This property follows
from the vanishing of certain integrals, see property
P5 in Sect.\ 7 of \cite{K}.
\end{enumerate}
The lowest order terms in $\epsilon$ are also the result
of explicit calculation of simplest Feynman integrals:
\begin{enumerate}
\item[C.] $P(\alpha)(f\otimes g)=fg+\epsilon \alpha(df\otimes fg)
+O(\epsilon^2)$.
\item[D.] $A(\xi,\alpha)=\xi+O(\epsilon)$.
\item[E.] $F(\xi,\eta)=O(\epsilon)$.
\end{enumerate}

It is worth noticing that the equations (i)--(iv) may be
written succinctly as a Maurer--Cartan equation 
$\delta S+\frac12[S,S]_G=0$ for $S=P+A+F$, where $[\,\cdot,\cdot]_G$
is the Gerstenhaber bracket on the Hochschild complex.

\section{Formal geometry}\label{s-fg}
In order to apply Kontsevich's formula to the jet bundle 
$E_0$ we need to introduce coordinates. This is done
using ideas of formal geometry, see \cite{K}: let 
 $\Mco$ be the manifold of 
jets of coordinates systems
on $M$. A point in $\Mco$ is an
 infinite jet  at zero of local diffeomorphisms
 $U\subset 
\mathbbm{R}^d\to M$ defined on some open neighborhood $U$ of $0\in\mathbbm{R}^d$. 
Two such maps define the same infinite
jet if and only if
 their Taylor expansions at zero (for any choice of local coordinates
on $M$) coincide. We have a projection $\pi:\Mco\to M$ 
sending $\varphi$
to $\varphi(0)$.

The group $\mathrm{GL}(d,\mathbbm{R})$ of linear diffeomorphisms
acts on $\Mco$ and we set $\Ma=\Mco/\mathrm{GL}(d,\mathbbm{R})$. Define
$\tilde
E_0
:=\Mco\times_{\mathrm{GL}(d,\mathbbm{R})}\mathbbm{R}[[y^{1},\dots,y^{d}]]$.
Moreover the Lie algebra $W_d$ of formal vector fields acts
on $\Mco$ by infinitesimal coordinate transformations. The action
is given by an isomorphism from $W_d$ to the tangent space at any
point of $\Mco$. The inverse map defines the 
{\em Maurer--Cartan form}
$\omega_{\mathrm{MC}}$, a $W_d$-valued one-form on $\Mco$.
The fact that this one-form comes from a Lie algebra action
implies that it obeys 
the Maurer--Cartan equation $d\omega_{\mathrm{MC}}+\frac12
[\omega_\mathrm{MC},\omega_\mathrm{MC}]=0$.
As a consequence, $\tilde D_0=d+\omega_\mathrm{MC}$ is a flat
connection on the trivial bundle $\Mco\times \cA$
over $\Mco$ and it is
easy to check that it descends to $\tilde E_0$.

We will need
the fact that the fibers of the bundle $\Ma\to M$ are contractible
so that there exist sections $\varphi^{\mathrm{aff}}:M\to \Ma$.
 For example,
the exponential map of a torsion free connection on the
tangent bundle defines such a section.

\begin{lemma}
 Let $\varphi^\mathrm{aff}:
x\mapsto [\varphi_x]$ be a section of $\Ma$.
The map sending a point $(x,f)$ of the jet bundle $E_0$
($f$ is the jet of a function at $x$) to the
class of
$(\varphi_x,\text{Taylor expansion of}\; f\circ\varphi_x)$ is
an isomorphism of vector bundles with flat connection
from  $(E_0,D_0)$ to $(\varphi^{\mathrm{aff}*}\tilde E_0,
\varphi^{\mathrm{aff}*}\tilde D_0)$.
\end{lemma}
{}From now on, we fix a section $\varphi^\mathrm{aff}$ and
identify $(E_0,D_0)$ with
  $(\varphi^{\mathrm{aff}*}\tilde E_0,
\varphi^{\mathrm{aff}*}\tilde D_0)$.
Here is an explicit local description of $E_0$ and $D_0$. On any open contractible
subset  $U\subset M$ we may choose a section $\varphi$ of $M^\mathrm{coor}\to
M$ such that $\pi\circ\varphi=\varphi^\mathrm{aff}$ where $\pi:\Mco\to\Ma$
is the canonical projection. We call such a section a {\em local lift} of
$\varphi^\mathrm{aff}$. Local lifts on $U$ differ by 
a $\mathrm{GL}(d,\R)$-gauge transformation $\varphi\mapsto\varphi\circ g$,
$g:U\to \mathrm{GL}(d,\R)$. A local lift induces a local trivialization
of $E_0|_U\simeq U\times \cA$, so that the isomorphism
$C^\infty(M)\to H^0(E_0,D_0)$ sends $f$ to the
Taylor expansion at zero of $f\circ\varphi$.
 One way to choose a local lift
is induced by a choice of local coordinates: let
$x^1,\dots,x^d:U\to \R$ be coordinates on $U\subset M$. 
Then there is a unique local lift $x\mapsto\varphi_x$
so that $\varphi^j_x=x^j\circ\varphi_x$
has the form
\[
\varphi^j_x(y)=x^j+y^j+\text{higher order terms in $y$}.
\]
The flat connection
$D_{0}$ on a local section $x\mapsto f_x(y)\in\R[[y^1,\dots,y^d]]$
is then
\[
D_0f_x(y)=\sum_{j=1}^ddx^j\left(\del{f_x(y)}{x^j}-\sum_{k,l=1}^d
T^l_j(x;y)\del{\varphi_x^k(y)}{x^l}
\del{f_x(y)}{y^k}
\right),
\]
where $T(x;y)$ is the matrix inverse to $(\partial\varphi^j_x(y)/\partial
y^k)$.
Indeed, it follows from the chain rule that Taylor
expansions of globally defined functions are
$D_0$-closed sections. 
Conversely, observe that 
$\sum_lT_{i}^{l}(x;y)\del{\varphi_x^k}{x^l}$ 
is a formal power series in $y$
beginning with
$\delta_{i}^{k}$ and whose coefficients are smooth in $x$. By this
property it follows immediately that the coefficients of a section
$\sigma$ of $E_0$ satisfying $D_{0}\sigma=0$ are determined by the
zeroth coefficient $\sigma^{0}(x)$. If we set
$\tilde\sigma=\sigma^{0}\circ\varphi$, we have
$D_{0}(\sigma-\tilde\sigma)=0$ and
$\left.(\sigma-\tilde\sigma)\right|_{y=0}=0$; but this implies
$\sigma=\tilde\sigma$. This shows that
 a section of $E_0$ is the Taylor
expansion of a globally defined function if
and only if it is $D_{0}$-closed.

\section{A deformation of the canonical connection}\label{s-4}

Let $E=E_0[[\epsilon]]$ be the bundle of jets of
$\R[[\epsilon]]$-valued functions. By definition, sections
of $E$ are formal power series in $\epsilon$ whose coefficients
are sections of $E_0$. We fix as above a section $\varphi^\mathrm{aff}$ of $\Ma$. Then $E$, with fiberwise Kontsevich product
associated to the Taylor expansion of $\alpha$, is a bundle
of $\R[[\epsilon]]$-algebras: choose a local lift 
$x\mapsto\varphi_x$ of $\varphi^\mathrm{aff}$ to a 
section of $\Mco$. This
induces a local trivialization of $E_0$ so that local
sections are given by $\cA[[\epsilon]]$-valued 
functions $x\mapsto f_x(y)$.
The product of local
sections $f,g$ is then the section
\[
x\mapsto(f\star g)_x=P(\alpha_x)(f_x\otimes g_x),\qquad
\alpha_x=\text{Taylor expansion of}\;
(\varphi_{x}^{-1})_*\alpha.
\]
It does not depend of the choice of local lift since $P$ is
$\mathrm{GL}(d,\R)$-equivariant, see Property A in Sect.~\ref{s-2}.

We now use $A$ of equation \eqref{e-A} to construct a deformation
$D=D_0+\epsilon D_1+\cdots$ of the canonical connection $D_0$.
Let $x\in M$ and $\xi\in T_xM$. Again we choose a local
lift $\varphi$ of $\varphi^\mathrm{aff}$ and set
\[\hat\xi_x=\varphi_x^*\omega_\mathrm{MC}(\xi)\in W_d.\]
The connection is then defined as
\[
D f_x=d f_x+A^M_xf_x,
\]
where $A_x^M(\xi)=A(\hat\xi_x,\alpha_x)$, $\xi\in T_xM$.
By properties A and B, this formula defines a connection on $E$
which is independent of the choice of local lift. By property
D, it is a deformation of $D_0$.

Equation (ii) in the formality theorem
 now implies the crucial statement:

\begin{proposition} For any
$f,g\in\Gamma(M,E)$, we have the Leibniz rule
\[D(f\star g)=Df\star g+f\star Dg.\]
\end{proposition}  
However $D$ is not flat, so that 
the algebra of horizontal sections with
respect to $D$ will not be isomorphic to $C^{\infty}(M)[[\epsilon]]$ as an
$\mathbbm{R}[[\epsilon]]$-module. We are going to discuss this
difficulty in the next section.

\section{Flattening the connection $D$}\label{s-5}

The connection $D$ can be extended to a 
(graded) derivation \[ D:\Omega^\Cdot(M,E)
\to\Omega^{\Cdot\,+1}(M,E)\] of the algebra
of differential forms on $M$ with values in $E$ (the
product on this algebra is defined by
the star-product in the fibers
and the wedge product on differential forms).
Its curvature $D^2$ is then an $\mathrm{End}(E)$ valued
two-form. Let $F^M\in\Omega^2(M,E)$ be the two-form
$x\mapsto F^M_x$ with
\[
F^M_x(\xi,\eta)=F(\hat\xi_x,\hat\eta_x,\alpha_x), \qquad\xi,\eta\in
T_xM.
\]
By properties A and B, 
$F^M$ is independent of the local lift
$\varphi_x$ needed to define $\hat\xi_x$,
$\hat\eta_x$ and $\alpha_x$.

Then the formality identities (iii), (iv) translate into
the following statements.

\begin{proposition}\label{p-34}
For any
$f\in\Gamma(M,E)$,
\begin{eqnarray}
D^2f&=&F^M\star f-f\star F^M,
\\
DF^M&=&0.\label{e-Bianchi}
\end{eqnarray}
\end{proposition}

The two-form
$F^M$ is called the \textit{Weyl curvature} of $D$. 
In general, a connection on a bundle of associative algebras with
the above properties, i.e. to be a derivation whose curvature
is an inner derivation such that its Weyl curvature satisfies
the Bianchi identity \eqref{e-Bianchi}, is called a \textit{Fedosov connection}.

We want now to modify $D$ so that it becomes flat still remaining
a derivation. The first observation is that
\[
\bar{D}:=D+[\gamma,\cdot]_{\star}
\]
is still a derivation for any
 $\gamma\in\Omega^1(M,E)$. The star-commutator is defined
by $[a,b]_\star=a\star b-b\star a$. Moreover, $\bar{D}$ turns out to
be again a Fedosov connection with Weyl curvature
$
\bar F^M=F^M+D\,\gamma+\gamma\star\gamma.
$
If we are able to find $\gamma$ so that $\bar F^M=0$, 
\begin{equation}\label{e-gamma}
F^M+D\,\gamma+\gamma\star\gamma=0,
\end{equation}
then $\bar{D}$-closed sections will form
a nontrivial subalgebra of $\Gamma(M,E)$.

The one-form $\gamma$ can be found as follows.
Since $F^M$ starts
at order $\epsilon$, see property E, we may write
$F=\epsilon F_1+\epsilon^2 F_2+\cdots$. The correction $\gamma$
does not need to have a term of order zero, so we write
$\gamma=\epsilon\gamma_1+\epsilon^2\gamma_2+\cdots$.
The equation $\bar F=0$ at order $\epsilon$ reads
\[
F_1+D_0\gamma_1=0,
\]
while the Bianchi identity imply
at this order $D_0F_1=0$. Hence, $D_0\gamma_1$ is equal to a $D_0$-closed
expression. 
But the $D_0$-cohomology is trivial in degree 2,
so it is possible to find a $\gamma_1$
that solves the equation. 
At higher order in $\epsilon$, one proves by induction that one always
has an equation of the form $D_0\gamma_k$ 
equal to a known $D_0$-closed form depending
on the lower order coefficients of $\gamma$ and $F^M$,
so a solution exists by the same argument.

Thus we have a flat connection $\bar D$ on the bundle
of algebras $E$, obeying the Leibniz rule. Let
$H^0(E,\bar D)$ be the $\R[[\epsilon]]$-algebra of
horizontal sections. We are left to show the
existence of a quantization map, an $\R$-linear homomorphism
\[
\rho:C^\infty(M)\simeq H^0(E_0,D_0)\to H^0(E,D),
\]
inducing an isomorphism 
$C^\infty(M)[[\epsilon]]\to H^0(E,\bar D)$ 
of $\R[[\epsilon]]$-modules. Then the
pull-back of the product defines a star-product on
$C^\infty(M)[[\epsilon]]$. To construct $\rho$ one
looks for a bundle map, also denoted by $\rho$, from
$E_0$ to $E$ of the form $\rho=\mathrm{Id}+\epsilon\rho_1+\dots$
 so that
\begin{equation}\label{e-rho}
\rho\circ D_0=\bar D\circ\rho,
\end{equation}
and that $\rho_i$ are given by differential operators.
Again, there are no cohomological obstructions to find such a
$\rho$: the $\rho_i$ can be found recursively by solving equations
of the form $D_0(\rho_i):=\rho_i\circ D_0-D_0\circ\rho_i=$ known.
The known term on the right is $D_0$-closed by the induction 
hypothesis, and $\rho_i$ can be found as a differential operator,
as the $D_0$-cohomology of the subbundle of $\mathrm{End}(E_0)$ given
by differential operators is trivial in degree 1.
Furthermore
 $\rho$ is unique if we require that $\rho(f)_x(y=0)=f_x(y=0)$,
for any section $f$.
The result is then:

\begin{thm} Let $(M,\alpha)$ be a Poisson manifold. Fix a
section $\varphi^\mathrm{aff}:M\to M^\mathrm{aff}$. Let $\star$,
$D$, $F^M$ be the corresponding product on $\Gamma(M,E)$, the connection
on $E$ and its curvature, respectively.
Let $\gamma\in\Omega^1(M,E)$, and $\rho\in \Gamma(M,\mathrm{Hom}(E_0,E))$ 
be solutions of \eqref{e-gamma},
\eqref{e-rho}, respectively.
Then
$\bar D=D+[\gamma,\cdot]_\star$ is a flat connection on $E$ obeying
the Leibniz rule $\bar D(f\star g)=\bar D(f)\star g
+f\star\bar D(g)$ and $\rho$
induces an isomorphism of $\R[[\epsilon]]$-modules $C^\infty(M)[[\epsilon]]
\to H^0(E,\bar D)$. The pull-back of the
product to $C^\infty(M)[[\epsilon]]$ is a star-product with first order term
$\alpha$.
\end{thm}

Thus the construction of a star-product for Poisson manifolds
requires solving an equation for $\gamma$ and an equation for 
$\rho$. We want to show that these equations can be solved in
an explicit way. The equation for $\gamma=\sum\epsilon^j\gamma_j$
can be solved recursively: at the $j$th step of the recursion,
one has an equation of the form
\begin{equation}\label{e-MissBertram}
D_0\gamma_j=\beta_j.
\end{equation}
where $\beta_j\in\Omega^2(M,E_0)$ obeys $D_0\beta_j=0$ 
and is known.
The point is that $\gamma_j$ may be constructed by a purely local
calculation: namely, let us introduce local coordinates and a
local trivialization of $E_0$ as in
Sect.~\ref{s-fg}. 

It is useful to define the total degree of a form on $M$
taking values in sections of $E_0$ as the sum of the form degree and the
degree in $y$. 
Then we write
\[
D_{0}=-\delta+D_{0}',
\]
where
\[
\delta:=\sum_{i=1}^d d x^i\del{}{y^i}
\]
is the zero-degree part and $D_{0}'$ has positive degree.
It follows immediately that $\delta^2=0$.
We can define a dual operator to $\delta$ on $E_0$-valued
differential forms: 
\begin{equation}\label{delta*}
\delta^*:=\sum_{i=1}^d y^i\,
\iota_{\frac\partial{\partial x^i}},
\end{equation}
where $\iota$ denotes inner multiplication.
It is easy to verify that $(\delta\delta^*+\delta^*\delta)\rho=k\rho$
for every form $\rho$ of total degree $k$. Thus, if we restrict to
$\delta$-closed 
forms of positive total degree $k$, we may then invert $\delta$ by
$\delta^{-1}\rho=\frac 1k\delta^*\rho$. This inverse yields the unique
form $\sigma$ such that $\delta\sigma=\rho$ and $\delta^*\sigma=0$.
This proves that the $\delta$-cohomology is concentrated in degree zero,
i.e., functions on $M$ (independent of $y$).

Then the solution of eq.~\eqref{e-MissBertram} obeying $\delta^*\gamma_i=0$
is
\[
\gamma_i=-\sum_{n=0}^\infty
(\delta^{-1}D_0')^n\delta^{-1}\beta_i.
\]
 The infinite sum converges in 
the sense of formal power series since the $n$-th term is
of degree at least $n$ in $y$.

The equation for $\rho$ can be solved similarly.
\section{Casimir and central functions}\label{s-6}
We describe variants of our constructions which are
relevant for the treatment of the center and for
the comparison with Fedosov's construction in the
symplectic case.

\subsection{Central closed two-forms}
Let $\omega=\omega_0+\epsilon\omega_1+\cdots
\in\Omega^2(M,E)$ be such that $
D\omega=0$ and $[\omega,f]_\star=0$ 
$\forall f\in\Gamma(M,E)$. Then we may 
construct a more general flat connection $\bar D$
by replacing \eqref{e-gamma} by
\begin{equation}\label{e-gammaomega}
F^M+D\,\gamma+\gamma\star\gamma=\omega.
\end{equation}
Indeed, 
the Bianchi identity holds also for $F^M-\omega$ and
$\bar D^2=0$ since $\omega$ is central. We get
thus a family of products parametrized by $D$-closed
central two forms $\omega$.

\subsection{The second construction of a quantization map}
The second variation concerns the quantization map.
It is possible to make it compatible with the center
in the following sense: let $Z_0(M)=\{f\in C^\infty(M)\,|\,
\{f,\cdot\}=0\}$ be the algebra
of Casimir functions.
There exists a quantization map $\rho:H^0(E,D_0)\to
H_0(E,\bar D)$ whose restriction
to $Z_0(M)$ induces an algebra isomorphism from
$Z_0(M)[[\epsilon]]$ 
onto the center of $H^0(E,\bar D)$. It is
constructed using the two remaining special cases
of the formality theorem involving the Poisson
vector field, vector fields and functions, see \cite{CFT}.

This quantization map can be used to construct a
linear map from the space of $D_0$-closed 
Casimir two-forms $B=\{\omega\in\Omega^2(M,E_0)\,|\,
D_0\omega=0, \{\omega,\cdot\}=0\}$ onto the
space of $\bar D$-closed central two-forms.
Thus we may parametrize our products by two-forms
in $B[[\epsilon]]$.

\subsection{Dependence on choices} If we
choose the homotopy $\delta^{-1}$ 
to determine canonical choices of $\gamma$ and $\rho$,
our construction of star-products 
depends on a Poisson bivector field $\alpha$,
a section $\varphi^\mathrm{aff}$ and a central closed
2-form in $B[[\epsilon]]$. 
It is possible to show that different
choices of $\varphi^\mathrm{aff}$ lead to equivalent
products (two star-products $\star$, $\star'$
are called equivalent if there is a 
series $\psi=\mathrm{Id}
+\epsilon\psi_1+\epsilon^2\psi_2+\cdots$ so that
$\psi(f\star g)=\psi(f)\star'\psi(g)$, $\forall
f,g\in C^\infty(M)$). Also, if we replace $\omega\in B[[\epsilon]]$
by $\omega+D_0\beta$ for a Casimir one-form $\beta\in\Omega^1(M,E_0)$, we obtain an equivalent star-product.
Still in general 
we do not get all star-products up to equivalence
in this way. To obtain all equivalence classes of star%
-products, we have to take $\alpha$ to be a power series
in $\epsilon$. Then \cite{K}
equivalence classes of star-products
are in one-to-one correspondence with $\R[[\epsilon]]$-%
valued Poisson vector fields $\alpha$ modulo
formal paths in the group of 
diffeomorphisms generated by flows of 
vector fields $\epsilon\xi_1+\epsilon^2\xi_2+\cdots$. 
This correspondence is realized by taking $\omega=0$
in our construction. If we take another $\omega$, we
get thus a product corresponding to $\omega=0$ 
for a different $\alpha$. This defines a map from
$B[[\epsilon]]$ modulo $D_0$-exact central two-forms to 
the equivalence classes of $\R[[\epsilon]]$-valued
Poisson bivector fields. It would be interesting to describe
this map. For infinitesimal $\omega$ and to leading order
in $\epsilon$ this map has the following description.
Let $C^{\Cdot,\Cdot}=
\Omega^\Cdot(M,\wedge^\Cdot \mathrm{Der}(E_0))$
be the double complex of differential forms with values
in the jets of formal multivector fields. The 
differentials are $D_0$ and the differential 
$\delta_\alpha$
of Poisson cohomology (the Schouten--Nijenhuis bracket
with $\alpha$) on  the fibers.
Then $\omega$ defines a class in $H^2(H^0(C,\delta_\alpha),D_0)$. As the $D_0$-cohomology is trivial except in 
degree zero, we have a map
\[
j:H^2(H^0(C^{\Cdot,\Cdot},\delta_\alpha),D_0)\to 
H^0(H^2(C^{\Cdot,\Cdot},\delta_\alpha),D_0)=H_\alpha^2(M),
\]
given by $\omega\mapsto \delta_\alpha D_0^{-1}\delta_\alpha D_0^{-1}\omega$.
Thus $j$ sends the class of $\omega$ to an element of
the second Poisson cohomology group, which consists of
equivalence classes of infinitesimal variations
of the Poisson bivector field $\alpha$.

\section{Examples, related results}\label{s-7}

\subsection{Open subsets of $\R^d$} Let us check that our construction
gives back the original Kontsevich formula for open subsets of $\R^d$.
In this case we may take $\varphi^\mathrm{aff}$ to be given by
\[
\varphi^j_x(y)=x^j+y^j.
\]
Then $D_0=\sum_{j=1}^d dx^j(\del{}{x^j}-\del{}{y^j})$. Since $A(\xi,\alpha)=\xi$ 
for constant $\xi$, we have $D=D_0$ and $F=0$. We may then choose
$\gamma=0$ and $\rho=\mathrm{Id}$ and the result is Kontsevich's
formula $f\star g=\sum_{n=0}^\infty\frac{ \epsilon^n}{n!}
U_n(\alpha,\dots,\alpha)(f\otimes g)$.

\subsection{Leaves of Poisson foliations \cite{ADGM}}
If $N$ is a foliated Poisson manifold so that the Poisson 
bivector field
is tangent to the leaves, then each leaf is a Poisson submanifold.
 An {\em adapted chart} with domain $V\subset M$ is a chart
$\psi:V\to \R^d$ so that the leaves
are given by equations $\psi^j(x)=\mathrm{const}$, $j=m+1,\dots,d$.
A star-product is called {\em tangential}
if on the domain $V$ of any adapted chart and for any $f\in C^\infty(M)$
so that $f|_V$ is constant on each leaf,
one has  $f\star g|_V=fg|_V$, $\forall g\in C^\infty(M)$.
In \cite{ADGM} it is shown that if $\varphi^\mathrm{aff}$ is chosen
to be {\em adapted} to the foliation, in the sense
that for all $x\in M$, there is
a representative  $\varphi_x$ of $\varphi^\mathrm{aff}(x)$
so that $\varphi_x^{-1}$ is an adapted chart,
then $\star$, with $\rho$ from the ``second construction'', is tangential.
Using this property, one can restrict
 star-products to the leaves of a foliated Poisson manifold. 
\subsection{The case of a symplectic manifold}
Suppose $(M,\alpha)$ is a Poisson manifold whose Poisson
bracket comes from a symplectic form $\Omega$.

A {\em Darboux section} of $\Ma$ is a
section $x\mapsto[\varphi_x]$ of $\Ma$ such that, for all
$x\in M$, $\varphi_x^*\Omega$ is a constant two-form.
Symplectic manifolds always have Darboux sections. For
example, a torsion free symplectic connection gives rise
to a Darboux section (but not through the
exponential map), see \cite{F}, Section 2.5.

In Darboux coordinates the Kontsevich product reduces to the
Moyal product, and also the other objects $A, F$ may be described 
explicitly. As a result, we have explicit formulae in terms
of a local lift of $\varphi^\mathrm{aff}$ with the property that
$\alpha_x=(\varphi_x^{-1})_*\alpha$ is a constant bivector
field $\alpha_0$.
In this case,
 $\varphi_x^*\omega_{MC}$ is a 1-form on $M$ with values in the
Hamiltonian formal vector fields. The corresponding hamiltonian
functions $h_x$, normalized by $h_x(y=0)=0$, define an $E_0$-valued
one-form $x\mapsto h_x$.
Then we have
\begin{enumerate}
\item[(i)]
$(f\star g)_x(y)=f_x\star g_x(y)=\exp\left(\epsilon\sum\alpha_0^{ij}
\frac\partial{\partial y_1^i}\frac\partial{\partial y_2^j}
\right)
f_x(y_1)g_x(y_2)|_{y_1=y_2=y}$.
\item[(ii)]
$Df_x=d_xf_x+\frac1{2\epsilon}
[h_x,f_x]_\star$.
\item[(iii)]
$F_x^M(\xi,\eta)=\frac1{4\epsilon^2}\left(
[\langle h_x,\xi\rangle,\langle h_x,\eta\rangle]_\star-
2\epsilon\{
\langle h_x,\xi\rangle,
\langle h_x,\eta\rangle\}\right)$.
\end{enumerate}
Moreover, in the symplectic case, the space $B$ of closed
Casimir two-forms consists of closed two-forms in $\Omega^2(M)[[\epsilon]]\subset\Omega^2(M,E)$.
Using these formulae the comparison with Fedosov's original
construction \cite{F} becomes clear. The latter starts with
the definition of the Weyl bundle:
let $F(M)$ be the principal $\mathrm{GL}(d,\R)$-bundle
of frames (bases of tangent spaces). The Weyl bundle is
the associated vector bundle 
$W=F(M)\times_{\mathrm{GL}(d,\R)}\R[[y^1,\dots,y^d]]
[[\epsilon]]$. 
The fiber over $x$ is the space of (formal) functions on
the tangent space at $x$, a symplectic vector space. So
there is a Moyal product in each fiber, and $W$ is a bundle
of algebras. A symplectic connection $\nabla$ on the
tangent bundle induces a connection $D_\mathrm{F}$
on $W$ obeying the Leibniz rule on sections. Its
curvature is $D_\mathrm{F}^2=(2\epsilon)^{-1}[R,\cdot]$, where 
the $W$-valued 2-form
$R=-\frac12\Omega(\nabla^2y,y)$ is the quadratic form on $TM$
associated with the curvature of $\nabla$. For any
closed two-form of the
form $\Omega(\epsilon)=-(2\epsilon)^{-1}\Omega+
\Omega_0+\epsilon\Omega_1+\cdots
\in\Omega^2(M)[[\epsilon]]$, 
Fedosov shows that the equation $D_\mathrm{F}\gamma_\mathrm{F}+
\gamma_\mathrm{F}\star\gamma_\mathrm{F}+R=\Omega(\epsilon)$  has a solution
$\gamma_\mathrm{F}$ of the form $(2\epsilon)^{-1}\Omega_{ij}y^idx^j+
\gamma_{\mathrm{F}1}+\epsilon\gamma_{\mathrm{F}2}+\cdots$,
obeying the normalization condition $\gamma_{\mathrm{F}}|_{y=0}=0$. It follows that connection $\bar D_\mathrm{F}=D_\mathrm{F}+[\gamma_\mathrm{F},\cdot]$ is flat.

Let $\nabla$ be the symplectic connection on $TM$ whose
Christoffel symbols are
obtained from the 2-jet of the Darboux section via
\[
\varphi^j_x(y)=x^j+y^j-\frac12\sum_{k,l} \Gamma^j_{kl}(x)y^ky^l+\cdots.
\] 
We have a bundle isomorphism $E\to W$ sending 
the jet at $x$ of a function $f$ to the class of
$(e,f\circ\varphi_x)$, where $\varphi_x$ is the unique
representative of $\varphi^{\mathrm{aff}}(x)$ that maps
the standard frame of $\R^d$ to the frame $e$. By (i)
this is an isomorphism of algebra bundles.  The connection
$\bar D_\mathrm{F}$ is then of the form $\bar D=
D+[\gamma,\cdot]$ with $\gamma_\mathrm{F}=
(2\epsilon)^{-1}h_x+\gamma$, where $\gamma$ is a solution
of \eqref{e-gammaomega} with $\omega=\Omega_0+\epsilon\Omega_1+\cdots$.

The result is that our star-product constructed using a
closed two-form $\omega\in\Omega^2(M)[[\epsilon]]$ 
is equivalent to Fedosov's star-product
associated to the class of 
$\Omega(\epsilon)=-(2\epsilon)^{-1}\Omega+\omega$. Details will be
presented elsewhere \cite{CFT2}.

\end{document}